\documentclass[12pt]{amsart}

\usepackage{amsmath}
\usepackage{latexsym}
\usepackage{amssymb}
\usepackage{amsfonts}
\usepackage{color}

\usepackage{graphics}

\renewcommand{\proof}{\par\noindent{\it Proof.\ \ }}
\def\qed{\ifmmode\square\else\nolinebreak\hfill
$\square$\fi\par\vskip12pt}

\renewcommand{\proof}{\par\noindent{\it Proof.\ \ }}
\def\qed{\ifmmode\square\else\nolinebreak\hfill
$\square$\fi\par\vskip12pt}

\def\ov{\overline} \def\un{\underline}
\def\l{\langle} \def\r{\rangle}

 \def\ZZ{\mathbb Z}

\def\S{{\rm S}} 
\def\J{{\sf J}} \def\M{{\sf M}} \def\soc{{\sf soc}}
  
 \def\mod{{\sf mod~}}
 
\def\val{{\sf val}}  
 \def\Aut{{\sf Aut}}

\def\Out{{\sf Out}}  \def\K{{\sf K}}

\def\Ga{{\it \Gamma}}

\def\Del{{\it\Delta}}
\def\Ome{{\it\Omega}}

\def\a{\alpha} \def\b{\beta}  \def\s{\sigma}
\def\t{\tau} 

\def\POmega{{\rm P\Omega}}
\def\Sp{{\rm Sp}} \def\PSp{{\rm PSp}}

\def\A{{\rm A}}

\def\PSL{{\rm PSL}}
  
\def\AGL{{\rm AGL}}

 \def\PSU{{\rm PSU}}

\newtheorem{theorem}{Theorem}[section]
\newtheorem{lemma}[theorem]{Lemma}%
\newtheorem{proposition}[theorem]{Proposition}%
\newtheorem{hypothesis}[theorem]{Hypothesis}%

    \def\ZZ{\mathbb Z}

\begin{document}

\title[$s$-arc-transitive digraphs]
{Vertex-primitive $s$-arc-transitive digraphs of alternating and symmetric groups}
\thanks{1991 MR Subject Classification 20B15, 20B30, 05C25.}
\thanks{This work was supported by
National Natural Science Foundation of China (11461007, 11231008).}

\author[Pan]{Jiangmin Pan}
\address{J. M. Pan\\
School of Statistics and Mathematics\\
Yunnan University of Finance and Economics\\
Kunming \\
P. R. China}
\email{jmpan@ynu.edu.cn}

\author[Wu]{Cixuan Wu}
\address{C. C. Wu\\
School of Statistics and Mathematics\\
Yunnan University of Finance and Economics\\
Kunming \\
P. R. China}

\author[Yin]{Fugang Yin}
\address{F. G. Yin\\
School of Statistics and Mathematics\\
Yunnan University of Finance and Economics\\
Kunming \\
P. R. China}

\maketitle

\begin{abstract}
A fascinating problem on digraphs is the existence problem of
the finite upper bound on $s$ for all vertex-primitive
$s$-arc-transitive digraphs except directed cycles (which is known to be
reduced to the almost simple groups case).
In this paper, we prove that $s\le 2$ for all $G$-vertex-primitive $s$-arc-transitive digraphs
with $G$ an (insoluble) alternating or symmetric group,
which makes an important progress towards a solution of the problem.
The proofs involves some methods that may be
used to investigate other almost simple groups cases.
\end{abstract}

\qquad {\textsc k}{\scriptsize \textsc {eywords.}} {\footnotesize
upper bound, $s$-arc-transitive digraph, vertex-primitive, alternating group}

\section{Introduction}

A {\it digraph} (directed graph) $\Ga$ is a pair $(V\Ga,\to)$
with vertex set $V\Ga$ and an antisymmetric irreflexive binary relation
$\to$ on $V\Ga$. All digraphs considered in this paper are finite.
For a positive integer $s$,
an $s$-arc of $\Ga$ is a sequence $v_0,v_1,\dots,v_s$
of vertices such that $v_i\to v_{i+1}$ for each $i=0,1,\dots,s-1$.
A $1$-arc is also simply called an {\it arc}.
A transitive permutation group $G$ is {\it primitive} on a set $\Ome$
if $G$ preserves no nontrivial partition of $\Ome$
(or equivalently, the point stabilizer
of $G$ is maximal in $G$).
For an automorphism group $G$ of $\Ga$,
we call that $\Ga$ is {\it $(G,s)$-arc-transitive}
if $G$ is transitive on the set of $s$-arcs of $\Ga$,
and $\Ga$ is {\it $G$-vertex-primitive} if $G$ is primitive on the vertex set of $\Ga$.
It is easy to see that $s$-arc-transitive digraphs with $s\ge 2$
are necessarily $(s-1)$-arc-transitive.

In sharp contrast with the undirected graphs,
where a well known result of Weiss \cite{Weiss}
states that finite undirected graphs other than cycles
can only be $s$-arc-transitive for $s\le 7$,
Praeger \cite{Praeger89} proved that there are infinite many
$s$-arc-transitive digraphs for unbounded $s$ other than directed cycles.
This interesting gap stimulated a series of constructions
\cite{CPW93,CLP95,Evans97,MMSZ02} for such digraphs
(which are called {\it highly transitive digraphs} in the literature).
However, although various constructions for $s$-arc-transitive digraphs
are known, finding vertex-primitive $s$-arc-transitive digraphs with $s\ge 2$
seems to be a very intractable problem:
in a survey paper of Praeger \cite{Praeger90} in 1990,
she said ``no such examples
have yet been found despite considerable effort by several people".
The existence problem of vertex-primitive 2-arc-transitive
digraphs besides directed cycles has just been solved until 2017 by Giudici, Li and Xia \cite{GLX18}
by constructing an infinite family of such digraphs with valency 6,
and no vertex-primitive $3$-arc-transitive digraphs have been founded yet.
These naturally motivate the following interesting problem
(posted by Giudici and Xia \cite{MX18}).

\vskip0.1in
\noindent{\bf Question A.} Is there an upper bound on $s$
for all vertex-primitive $s$-arc-transitive digraphs that are not directed cycles?
\vskip0.1in

A group $G$ is said to be {\it almost simple}
if there is a nonabelian simple group $T$ such that $T\lhd G\le\Aut(T)$.
A systematic investigation of the O'Nan-Scott types of primitive permutation groups
has reduced Question A to the almost simple case
by proving that an upper bound on $s$ for vertex-primitive
$s$-arc-transitive digraphs $\Ga$ with $\Aut\Ga$ almost simple will be an upper
bound on $s$ for all vertex-primitive $s$-arc-transitive digraphs,
see \cite[Corollary 1.6]{MX18}.
Thus a reasonable strategy for Question A is to
investigate the upper bound of $s$ for all simple groups (the sporadic simple groups case
can generally be done especially with the help of the computer program).
In this paper, we will do such for
alternating and symmetric groups.
Our main result is as following.

\begin{theorem}\label{Thm-1}
Let $\Ga$ be a $G$-vertex-primitive $(G,s)$-arc-transitive digraph,
where $G$ is an insoluble alternating or symmetric group.
Then $s\le 2$.
\end{theorem}

In light of Theorem~\ref{Thm-1} and \cite[Theorem 1.2]{GLX17} for linear groups,
we would like to propose the following conjecture.

\vskip0.1in
\noindent{\bf Conjecture.} The upper
bound on $s$ for all vertex-primitive $s$-arc-transitive digraphs
except direct cycles is $2$.

\vskip0.1in
The layout of this paper is as follows.
We give some preliminary results in Section 2.
Notice that, for the digraphs in Theorem~\ref{Thm-1},
the vertex stabilizers of $G$ satisfy parts (a)-(f) of Theorem~\ref{Max-SubG} below
(obtained by Liebeck, Praeger and Saxl).
Parts (a) and (c) are investigated in Section 3,
and part (d) is considered in Section 4.
In the final Section 5,
we analyse parts (b,e,f) and complete the proof of Theorem~\ref{Thm-1}.

\section{Preliminaries}

Throughout the paper, we always use the following notations,
where $G$ is a group, $n$ is a positive integer and $p$ is a prime.

\vskip0.1in

$\pi(G)$: the set of prime divisors of the order of $G$.
\vskip0.05in
$\pi(n)$: the set of prime divisors of $n$.

\vskip0.05in
$n_p$: the maximal power of $p$ dividing $n$.
\vskip0.05in
$\soc(G)$: the socle of $G$, namely the product of all minimal normal subgroups of $G$.


\vskip0.1in
The following result is a consequence (also easy to prove directly)
of the so-called Legendre's formula,
that will be used repeatedly in this paper.

\begin{lemma}\label{p-part}
For each positive integer $n$ and prime $p$,
we have $(n!)_p<p^{n\over p-1}$.
\end{lemma}

For positive integers $a\ge 2$ and $m\ge 2$, a prime $r$ is called
a {\it primitive prime divisor} of $a^m-1$
if $r$ divides $a^m-1$ but not divides $a^i-1$ for each $i=1,2,\dots,m-1$.
The next is a well-known theorem of Zsigmondy,
see \cite[Theorem IX.8.3]{BB82}, where the last statement
follows easily by the Fermat's Little Theorem.

\begin{lemma}\label{Zsigmondy}
For positive integers $a\ge 2$ and $m\ge 2$,
$a^m-1$ has a primitive prime divisor $r$
if $(a,m)\ne (2,6)$ and $(2^e-1,2)$ with $e\ge 2$ an integer.
Moreover, $r\equiv 1(\mod m)$, in particular $r>m$.
\end{lemma}

The following theorem of Liebeck, Praeger and Saxl \cite{LPS87} determines the maximal subgroups of
alternating and symmetric groups which provides a starting point of this paper.

\begin{theorem}\label{Max-SubG}
Let $G=\A_n$ or $\S_n$,
and $H\ne\A_n$ a maximal subgroup of $G$.
Then $H$ satisfies one of the following:

\begin{itemize}
\item[(a)] $H=(\S_m\times \S_k)\cap G$, with $n=m+k$ and $m<k$ (intransitive case);
\item[(b)] $H=(\S_m\wr\S_k)\cap G$, with $n=mk$, $m>1$ and $k>1$ (imprimitive case);
\item[(c)] $H=\AGL(k,p)\cap G$, with $n=p^k$ and $p$ prime (affine case);
\item[(d)] $H=(T^k.(\Out(T)\times\S_k))\cap G$,
with $T$ a nonabelian simple group, $k\ge 2$ and $n=|T|^{k-1}$ (diagonal case);
\item[(e)] $H=(\S_m\wr\S_k)\cap G$, with $n=m^k$, $m\ge 5$ and $k>1$ (wreath case);
\item[(f)] $T\lhd H\le\Aut(T)$, with $T$ a nonabelian simple group,
$T\ne\A_n$ and $H$ acts primitively on $\Ome$ (almost simple case).
\end{itemize}
\end{theorem}

We remark that not all the groups $H$ satisfying one of parts (a-f)
are exactly maximal subgroups of $G$, namely there have a few exceptions,
see \cite[T{\tiny HEOREM}]{LPS87}.
The next known result presents a necessary and sufficient condition
of $s$-arc-transitivity of digraphs, refer to \cite[Lemma 2.2]{MX18}.

\begin{lemma}\label{iff}
Let $\Ga$ be a digraph,
and $v_0\to v_1\to \cdots\to v_{s-1}\to v_s$ be an $s$-arc of $\Ga$ with $s\ge 2$.
Suppose $G\le\Aut\Ga$ acts arc-transitively on $\Ga$.
Then $G$ acts $s$-arc-transitively on $\Ga$ if and only if
$$G_{v_1v_2\dots v_i}=G_{v_0v_1\dots v_i}G_{v_1\dots v_iv_{i+1}},~
for~ each~i\in\{1,2,\dots,s-1\}.$$
\end{lemma}

\vskip0.1in
For a group $G$, an expression $G=HK$ with $H$ and $K$ being subgroups of $G$
is called a {\it factorization} of $G$, where $H$ and $K$ are called {\it factors} of $G$.
In particular, $G=HK$ is called
a {\it homogeneous factorization}
if $H$ is isomorphic to $K$,
and is called a {\it maximal factorization} if both $H$ and $K$ are maximal subgroups of $G$.


\begin{lemma}\label{AS-Factori}{\rm(\cite[Proposition 3.3]{GLX17})}
Let $G$ be an almost simple group with socle $T$.
Suppose $G=AB$ is a homogeneous factorization. Then one of the following holds.
\begin{itemize}
\item[(a)] Both $A$ and $B$ contain $T$.
\item[(b)] $A$ and $B$ are almost simple groups with socles both isomorphic to $S$, where $(T,S)$
lies in the following table, where $q$ is a prime power and $f>2$.
\end{itemize}
\end{lemma}

\begin{table}[ht]
\[\begin{array}{lllll} \hline
T & \A_6 & \M_{12} &  \Sp_4(2^f) & \POmega_8^+(q)  \\ \hline
S & \A_5 & \M_{11} & \Sp_2(4^f) & \Omega_7(q)   \\ \hline
\end{array}\]
\end{table}

\begin{lemma}\label{Wreath-Factor}{\rm(\cite[Lemma 3.5]{GLX17})}
Let $R\wr\S_k$ be a wreath product with base group $M=R^k=R_1\times\cdots\times R_k$,
and $T\wr\S_k\le G\le R\wr\S_k$ with $T\le R$.
Suppose $G=AB$ is a homogeneous factorization of $G$
such that $A$ is transitive on $\{R_1,\dots,R_k\}$.
Denote by $\phi_i(A\cap M)$ the projection of $A\cap M$ on $R_i$
for $i=1,2,\dots,k$.
Then $\phi_1(A\cap M)=\cdots=\phi_k(A\cap M)$
and $\pi(T)\subseteq\pi(\phi_1(A\cap M))$.
\end{lemma}

\begin{lemma}\label{Factor-1}
Let $G$ be an almost simple group with socle $T=\PSL_k(q)$,
where $k\ge 2$ and $q=p^e$ is a prime power.
If $G=HK$ with $H$ and $K$ subgroups of $G$ such that
$\pi(H)\cap\pi(K)\supseteq \pi(G)\setminus\pi(p(p-1))$,
then either
\begin{itemize}
\item[(i)] at least one of $H$ and $K$ contains $T$; or
\item[(ii)] $k=2$ and $q=2^e-1\ge 7$ is a Mersenne prime.
\end{itemize}
\end{lemma}

\proof Let $H_0$ and $K_0$ be maximal subgroups
of $G$ containing $H$ and $K$, respectively.
Then $G=H_0K_0$ is a maximal factorization.
Such factorizations for $G$ being an almost simple group
with socle $\PSL(d,q)$ are classified in \cite[TABLE 1]{LPS90}.
By checking the list, the lemma follows.\qed

We give an observation to end this section. Denote by $\val(\Ga)$
the valency of a regular digraph $\Ga$.

\begin{lemma}\label{s-ArcT}
Let $\Ga$ be a $(G,s)$-arc-transitive digraph with $G\le\Aut\Ga$ and $s\ge 1$.
Then $\val(\Ga)^s\mid |G_v|$ for each $v\in V\Ga$.
\end{lemma}

\proof Set $m=\val(\Ga)$, and let $v=v_0\to v_1\to\cdots\to v_s$ be an $s$-arc of $\Ga$.
Since $\Ga$ is  $(G,s)$-arc-transitive,
$G_{v_0v_1\dots v_{i-1}}$ is transitive on the out-neighbor set $\Ga^+(v_{i-1}):=\{u\in V\Ga\mid v_{i-1}\to u\}$
for each $i=1,2,\dots,s$.
Then as $|\Ga^+(v)|=m$ for each $v\in V\Ga$,
we deduce
$|G_{v_0v_1\dots v_{i-1}}: G_{v_0v_1\dots v_{i}}|=m$.
It follows $|G_v|=|G_{v_0}|=m^s|G_{v_0v_1\dots v_{s}}|$,
the lemma follows.\qed

\section{Subgroups $(a)$ and $(c)$}

For convenience, we always suppose the following hypothesis holds in the rest of this paper.

\begin{hypothesis}\label{Hypo}
Let $\Ga$ be a $G$-vertex-primitive $(G,s)$-arc-transitive digraph
of valency at least three, where $s\ge 1$
and $G=\A_n$ or $\S_n$ with $n\ge 5$ is an automorphism group
of $\Ga$. Take an arc $u\to v$ of $\Ga$,
and let $g\in G$ such that $u^g=v$ and set
$w=v^g$. Then $u\to v\to w$ is a $2$-arc of $\Ga$.
Set $\Ome=\{1,2,\dots,n\}$.
Then $G$ acts naturally on $\Ome$.
\end{hypothesis}

Under Hypothesis~\ref{Hypo}, $G_{vw}=G_{uv}^g$
and $G_v$ is a maximal subgroup of $G$.
Hence $G_v$ satisfies parts $(a)-(f)$ of Theorem~\ref{Max-SubG}.
In this section, we investigate the cases where $G_v$ satisfies parts $(a)$
and $(c)$.

\begin{lemma}\label{part-a}
Suppose $G_v$ satisfies part $(a)$ of Theorem~$\ref{Max-SubG}$.
Then $s=1$.
\end{lemma}

\proof Suppose for a contradiction that $s\ge 2$.
By Theorem~\ref{Max-SubG}, $G_v\cong(\S_m\times \S_k)\cap G$, with $n=m+k$ and $m<k$.
If $m=1$, then $G$ is 2-transitive on $V\Ga$,
so $\Ga$ is an undirected complete graph,
a contradiction.

Thus assume $m\ge 2$ in the following.
Since $s\ge 2$, $G_v=G_{uv}G_{vw}$ by Lemma~\ref{iff}
and $G_{vw}=G_{uv}^g\cong G_{uv}$.
Notice that $G$ has unique conjugate class of $(\S_m\times \S_k)\cap G$,
the action of $G$ on $V\Ga$ is permutation equivalent to
the natural induced action of $G$ on $\Ome^{\{m\}}$,
the set of $m$-subsets of $\Ome$.
We may thus identify $V\Ga$ with $\Ome^{\{m\}}$
and set $v=\Del:=\{1,2,\dots,m\}$.
Then $G_{v}=(\S_{\{1,\dots,m\}}\times\S_{\{m+1,\dots,n\}})\cap G$.
Clearly, $G_v=\S_{\{1,\dots,m\}}\times\S_{\{m+1,\dots,n\}}$ if $G=\S_n$,
and $G_v=(\A_{\{1,\dots,m\}}\times\A_{\{m+1,\dots,n\}}):\ZZ_2$ if $G=\A_n$.

Assume first $\Del\cap \Del^g=\phi$.
Without loss of generality,
we may suppose $\Del^g=\{m+1,m+2,\dots,2m\}$.
Then
$$G_{w}=G_{v}^g=(\S_{\{m+1,\dots,2m\}}\times\S_{\{1,\dots,m,2m+1,\dots,n\}})\cap G,~and$$
$$G_{u}=G_{v}^{g^{-1}}=(\S_{\{j_1,\dots,j_m\}}\times\S_{\Ome\setminus\{j_1,\dots,j_m\}})\cap G,$$
where $\{j_1,\dots,j_m\}\subseteq\Ome\setminus\{1,\dots,m\}.$
It follows
$$G_{uv}=G_{u}\cap G_{v}=(\S_{\{1,\dots,m\}}\times\S_{\{j_1,\dots,j_m\}}\times\S_{\{m+1,\dots,n\}\setminus\{j_1,\dots,j_m\}})\cap G,~and$$
$$G_{vw}=G_{v}\cap G_{w}=(\S_{\{1,\dots,m\}}\times\S_{\{m+1,\dots,2m\}}\times\S_{\{2m+1,\dots,n\}})\cap G.$$

\noindent Since $G_v=G_{uv}G_{vw}$,
we conclude

\vskip0.1in
$\S_{\{m+1,\dots,n\}}\cap G$

$=((\S_{\{j_1,\dots,j_m\}}\times\S_{\{m+1,\dots,n\}\setminus\{j_1,\dots,j_m\}})\cap G)
((\S_{\{m+1,\dots,2m\}}\times\S_{\{2m+1,\dots,n\}})\cap G).\hspace{30pt}(1)$

\vskip0.1in
Since $n=m+k\ge 5$ and $m<k$, we have $n-m=k\ge 3$.
If $n-m=3$ or $4$, one easily verifies Equation (1) is impossible, a contradiction.
If $n-m\ge 5$, then $\S_{\{m+1,\dots,n\}}\cap G\cong\A_{n-m}$ or $\S_{n-m}$ is almost simple,
and the two factors in the right side of Equation (1)
are conjugate in $G$ and so isomorphic,
by Lemma~\ref{AS-Factori},
the only possibility is $n-m=6$ and $m=5$,
and so $n=11$. Then Equation (1) leads to
$$\S_{\{6,\dots,11\}}\cap G=(\S_{\{j_1,\dots,j_5\}}\cap G)(\S_{\{6,\cdots,10\}}\cap G),~with~\{j_1,\dots,j_5\}\subseteq\{6,\cdots,11\}.\hspace{15pt}(2)$$
Clearly, one of $j_1,\dots,j_5$ equals $11$.
Then the intersection of the two factors in the right side of Equation (2) is isomorphic to $\A_4$ or $\S_4$,
it follows that the order of the group in the right side of Equation (2)
in a multiple of $25$, but the order of the group in the left side is not,
 also a contradiction.

Now assume $\Del\cap \Del^g\ne\phi$.
We may assume $\Del\cap\Del^g=\{1,\dots,l\}$ with $l<m$.
A direct computation shows
$$G_{uv}=(\S_{\{h_1,\dots,h_l\}}\times\S_{\{1,\dots,m\}\setminus\{h_1,\dots,h_l\}}
\times\S_{\{k_{l+1},\dots,k_m\}}\times\S_{\{m+1,\dots,n\}\setminus \{k_{l+1},\dots,k_m\}})\cap G,$$
$$G_{vw}=(\S_{\{1,\dots,l\}}\times\S_{\{l+1,\dots,m\}}
\times\S_{\{j_{l+1},\dots,j_m\}}\times\S_{\{m+1,\dots,n\}\setminus \{j_{l+1},\dots,j_m\}})\cap G,$$
where $\{h_1,\dots,h_l\}\subseteq \{1,\dots,m\}$,
$\{k_{l+1},\dots,k_m\}$
and $\{j_{l+1},\dots,j_m\}$ are both subsets of $\{m+1,\dots,n\}$.
Because $G_{v}=G_{uv}G_{vw}$,
we conclude
$$\S_{\{1,\dots,m\}}\cap G=((\S_{\{1,\dots,l\}}\times\S_{\{l+1,\dots,m\}})\cap G)
((\S_{\{h_1,\dots,h_l\}}\times\S_{\{1,\dots,m\}\setminus\{h_1,\dots,h_l\}})\cap G).$$
If $2\le m\le 4$, a simple computation may draw a contradiction.
If $m\ge 5$, 
by Lemma~\ref{AS-Factori},
the only possibility is $m=6$ and $l=1$,
then similar arguments as above may draw a contradiction.\qed

\begin{lemma}\label{part-c}
Suppose $G_v$ satisfies part $(c)$ of Theorem~$\ref{Max-SubG}$.
Then $s=1$.
\end{lemma}

\proof Suppose for a contradiction that $s\ge 2$.
By assumption, $G_v\cong\AGL(k,p)\cap G$ with $n=p^k$ and $p$ a prime,
and so $\soc(G_v)\cong\ZZ_p^k$.
Since $n\ge 5$, $(k,p)\ne(2,2)$.
If $(k,p)=(2,3)$, then $G_v\cong\ZZ_3^2:2\A_4$ or $\ZZ_3^2:2\S_4$,
a direct computation by Magma \cite{Magma} shows that
$G_v$ has no homogeneous factorization $G_v=G_{uv}G_{vw}$
with $|G_v:G_{uv}|\ge 3$, a contradiction.

Thus assume in the following $(k,p)\ne (2,2)$ and $(2,3)$.
Then $G_v$ is insoluble.
Let $M$ be a normal subgroup of $\AGL(k,p)$
such that $M\cong\ZZ_p^k:\ZZ_{p-1}$.
Set $\ov{G_v}=G_vM/M$, $\ov{G_{uv}}=G_{uv}M/M$
and $\ov{G_{vw}}=G_{vw}M/M$.
Then $\ov{G_v}$ is almost simple with $\soc(\ov{G_v})\cong\PSL_k(p)$.
Since $G_v=G_{uv}G_{vw}$
and $G_{uv}\cong G_{vw}$,
we have $\pi(G_{uv})=\pi(G_{vw})=\pi(G_v)$
and $\ov{G_v}=\ov{G_{uv}}\hspace{1pt}\ov{G_{vw}}$.
It follows that
$$\pi(\overline{G_{uv}})\cap\pi(\overline{G_{vw}})
\supseteq(\pi(G_{uv})\cap\pi(G_{vw}))\setminus\pi(M)\supseteq{\pi(\overline{G_v})\setminus\pi(p(p-1))}.$$
By Lemma~\ref{Factor-1}, either
\vskip0.1in
\begin{itemize}
\item[(i)] at least one of $\ov{G_{uv}}$ and $\ov{G_{vw}}$ contains $\soc(\ov{G_v})$; or
\item[(ii)] $k=2$ and $p=2^e-1\ge 7$ is a Mersenne prime.
\end{itemize}
\vskip0.1in
For case (i), without loss of generality, we may suppose
$\ov{G_{uv}}\supseteq\soc({\ov{G_v}})$.
Since $\ov{G_v}$ is almost simple,
$\soc(\ov{G_{uv}})=\soc({\ov{G_v}})$.
Then as $G_{vw}\cong G_{uv}\cong (M\cap G_{uv}).\ov{G_{uv}}$
and $M$ is soluble, both $\ov{G_{uv}}$ and
$\ov{G_{vw}}$ have the same unique insoluble composition factor $\PSL_k(p)$.
Since
$\ov{G_{vw}}/(\ov{G_{vw}}\cap\soc({\ov{G_v}}))\cong \ov{G_{vw}}\soc({\ov{G_v}})/\soc({\ov{G_v}})\le \ov{G_v}/\soc({\ov{G_v}})$ is soluble,
$\PSL_k(p)$ is a composition factor of $\ov{G_{vw}}\cap\soc({\ov{G_v}})$,
so $\ov{G_{vw}}\cap\soc({\ov{G_v}})=\soc({\ov{G_v}})$,
and hence $\soc(\ov{G_{vw}})=\soc({\ov{G_v}})$
as ${\ov{G_v}}$ is almost simple.
As $G_v=G_{uv}G_{vw}$,
at least one of $G_{uv}$ and $G_{vw}$, say $G_{uv}$,
has nontrivial intersection with $\soc(G_v)$,
then since $\ov{G_{uv}}\supseteq\soc(\ov{G_v})\cong\PSL_k(p)$
which acts irreducibly on $\soc(G_v)\cong\ZZ_p^k$,
we further conclude $G_{uv}\supseteq \soc(G_v)$.
Consequently, one easily has $\soc(G_{vw})=\soc(G_v)$ as $|G_{vw}|=|G_{uv}|$.
It follows
$$(\soc(G_v))^g=(\soc(G_{uv}))^g=\soc(G_{uv}^g)=\soc(G_{vw})=\soc(G_v).$$
Notice that $\l G_v,g\r=G$, we have $\soc(G_v)\lhd G$,
that is, $G_v$ contains a nontrivial normal subgroup of $G$,
hence $G$ acts unfaithfully on $V\Ga$, a contradiction.

For case (ii), $\ov{G_v}\cong\PSL_2(p).o$ with $o=1$ or $\ZZ_2$.
If one of $\ov{G_{uv}}$ and $\ov{G_{vw}}$ contains $\soc(\ov{G_v})\cong\PSL_2(p)$,
the discussions in the previous paragraph may draw a contradiction.
Suppose now none of $\ov{G_{uv}}$ and $\ov{G_{vw}}$ contains $\soc(\ov{G_v})$.
Since $\ov{G_v}=\ov{G_{uv}}\hspace{1pt}\ov{G_{vw}}$,
we may suppose (interchange $\ov{G_{uv}}$ and $\ov{G_{vw}}$
if necessary) that $\ov{G_{uv}}\le\ZZ_p:\ZZ_{p-1\over 2}.o$,
hence $|\ov{G_{uv}}|_2\le |o|$ as $p=2^e-1$.
Since $\ov{G_{uv}}\cong G_{uv}/(G_{uv}\cap M)$,
we conclude $|G_{uv}|_2\le |o||M|_2=2|o|$.
It follows
$$2^{e+1}|o|=2|\PSL_2(2^e-1).o|_2=|G_v|_2\le |G_{uv}|_2^2\le 4|o|^2\le 8|o|,$$
implying $e\le 2$ and so $p\le 3$, which is a contradiction.\qed

\section{Subgroups $(d)$}

Suppose Hypothesis~\ref{Hypo} holds. In this section, we consider that case where
$G_v$ satisfies part (d) of Theorem~\ref{Max-SubG},
namely $G_v\cong(T^k.(\Out(T)\times\S_k))\cap G$,
with $T$ a nonabelian simple group, $k\ge 2$ and $n=|T|^{k-1}$.
The main result of this section is the following assertion.

\begin{lemma}\label{part(d)}
Suppose $G_v$ satisfies part $(d)$ of Theorem~$\ref{Max-SubG}$.
Then $s\le 2$.
\end{lemma}

Lemma~\ref{part(d)} will be proved
by a series of lemmas in which 
the following proposition of Liebeck, Praeger and Saxl
(see \cite[P. 296, Corollary 5]{LPS2000}) will be used.
We remark there are two minor problems in Table 10.7 there,
namely `$p\le c$' in Row 1 should be `$p<c$' (for if $p=c$ then $L\supseteq T$),
and `$\ZZ_p:\ZZ_{p-1\over 2}$' in Row 7 should be
`$\ZZ_p:\ZZ_{p-1}$'.

\begin{proposition}\label{Coro5}
Let $G$ be an almost simple group with socle $T$.
Suppose that $L$ is a subgroup of $G$ such that $\pi(T)\subseteq \pi(L)$.
Then either

\begin{itemize}
\item[(i)] $T\subseteq L$; or
\item[(ii)] the possibilities for $T$ and $L$ are given in Table $1$.
\end{itemize}
\end{proposition}

\begin{table}
\[\begin{array}{llll} \hline
Row & T &  L\cap T  &  Remark  \\ \hline
1& \A_c &  \A_l\lhd L\le\S_l\times\S_{c-l}  &  p~ prime, p<c\Rightarrow~ p\le l  \\
2& \A_6 &  \PSL(2,5) & \\
3& \PSp_{2m}(q)~(m,q~even) &  L\rhd \Ome_{2m}^{-}(q)  &  \\
4& \POmega_{2m+1}(q)~(m~even,q~odd) &  L\rhd \Ome_{2m}^{-}(q)  &  \\
5& \POmega_{2m}^+(q)~(m~even) &  L\rhd\Ome_{2m-1}(q) &   \\
6& \PSp_4(q) &  L\rhd\PSp_2(q^2)  &   \\
7& \PSL_2(p)~(p=2^m-1) &  L\le\ZZ_p:\ZZ_{p-1}  & G=T.2  \\
8& \PSL_2(8)&  7.2, P_1  & G=T.3 \\
9& \PSL_3(3)&  13:3  &  G=T.2 \\
10& \PSL_6(2) &  P_1, P_5, \PSL_5(2)  &   \\
11& \PSU_3(3) &  \PSL_2(7) &   \\
12& \PSU_3(5) &  \A_7 &   \\
13& \PSU_4(2) &  L\le 2^4.\A_5~ or~\S_6 &  \\
14& \PSU_4(3) &  \PSL_3(4),\A_7 &   \\
15& \PSU_5(2) &  \PSL_2(11) &   \\
16& \PSU_6(2) &  \M_{22} &   \\

17& \PSp_4(7) &  \A_7 &  \\

18& \Sp_4(8) & ^2B_2(8) & G=T.3  \\

19& \Sp_6(2) & \S_8,\A_8,\S_7,\A_7 & \\

20& \POmega_8^+(2) & L\le P_i~(i=1,3,4), \A_9 & \\
21& G_2(3) & \PSL_2(13) &   \\
22& ^2F_4(2)' & \PSL_2(25) &  \\
23& \M_{11} & \PSL_2(11) & \\
24& \M_{12} & \M_{11}, \PSL_2(11) & \\
25 & \M_{24} & \M_{23} & \\
26 & HS & \M_{22} & \\
27 & McL & \M_{22} & \\
28 & Co_2 & \M_{23} & \\
29 & Co_3 & \M_{23} &  \\
 \hline
\end{array}\]
\caption{Subgroups $L$ with $\pi(T)\subseteq\pi(L)$ and $T{\not\subseteq}L$}
\end{table}

Let $M$ be a normal subgroup of $G_v$ isomorphic to $T^k.\Out(T)$.
Set $\ov{G_v}=G_vM/M$,
$\ov{G_{uv}}=G_{uv}M/M$ and
$\ov{G_{vw}}=G_{vw}M/M$.
Then $\A_k\le\ov {G_v}\le\S_k$.

\begin{lemma}\label{part-d}
If $\Ga$ is $(G,2)$-arc-transitive, then
$($interchange $\ov{G_{uv}}$ and $\ov{G_{vw}}$ if necessary$)$
$\ov{G_{uv}}\le\S_k$ is transitive
and $\phi_1(G_{uv}\cap M)=\cdots=\phi_k(G_{uv}\cap M)$.
Further, either $\phi_i(G_{uv}\cap M)\supseteq T$,
or the couple $(T,\phi_i(G_{uv}\cap M))$ $($as $(T,L)$ there$)$ satisfies Table $1$ of Proposition~$\ref{Coro5}$,
where $\phi_i(G_{uv}\cap M)$ denotes the projection of $G_{uv}\cap M$ on the $i$-component of $M$.
\end{lemma}

\proof Since $\Ga$ is $(G,2)$-arc-transitive, $G_v=G_{uv}G_{vw}$,
and hence $\ov{G_v}={\ov{G_{uv}}}\hspace{1pt}{\ov{G_{vw}}}$.
Noting that $\ov{G_v}\cong\A_k$ or $\S_k$,
by \cite[Lemma 2.3]{GLX17},
at least one of $\ov{G_{uv}}$ or $\ov{G_{vw}}$, say $\ov{G_{uv}}$,
is a transitive subgroup of $\S_k$.
It then follows from Lemma~\ref{Wreath-Factor} that
$\phi_1(G_{uv}\cap M)=\cdots=\phi_k(G_{uv}\cap M)$,
and $\pi(T)\subseteq\pi(\phi_1(G_{uv}\cap M))$.
Now by Proposition~\ref{Coro5}, the lemma follows.\qed

The following lemma treats the case
where $\phi_1(G_{uv}\cap M)$ contains $T$.

\begin{lemma}\label{lem-1}
Assume $T\subseteq \phi_1(G_{uv}\cap M)$. Then $s\le 2$.
\end{lemma}

\proof Suppose on the contrary $s\ge 3$.
By Lemma~\ref{part-d}, we may assume that
$\phi_1(G_{uv}\cap M)=\dots=\phi_k(G_{uv}\cap M)=T.o$
with $o\le\Out(T)$, and $\ov{G_{uv}}\le\S_k$ is transitive.
It follows that $G_{uv}\cap M$ has a unique insoluble composition factor $T$
with multiplicity, say $l$, dividing $k$.

We first prove $l=k$. If not,
then $l\le {k\over 2}$ as $l\mid k$.
It is known (or see \cite[P. 297, Corollary 6]{LPS2000})
that there is a prime $r\ge 5$
such that $r\mid |T|$ but $r$ does not divide $|\Out(T)|$.
So $|G_{uv}\cap M|_r=|T|_r^l\le |T|_r^{k/2}$,
and hence
$$\val(\Ga)_r={|G_v|_r\over |G_{uv}|_r}\ge{|G_v|_r\over |G_{uv}\cap M|_r|\S_k|_r}
\ge {|T|_r^k(k!)_r\over |T|_r^{k/2}(k!)_r}=|T|_r^{k/2}.$$
Since $s\ge 3$,
by Lemma~\ref{s-ArcT}, $\val(\Ga)_r^3\le |G_v|_r$,
it follows $|T|_r^{3k/2}\le |T|_r^k(k!)_r$.
However, as $(k!)_r<r^{k\over r-1}$ by Lemma~\ref{p-part},
we conclude $|T|_r^{k/2}<r^{k\over r-1}$,
which is a contradiction as $r\ge 5$.

Thus $l=k$. Consequently, $G_{uv}\cap M$ contains $\soc(G_v)\cong T^k$.
Since $\ov{G_{uv}}$ is transitive,
we further conclude
$\soc(G_v)$ is the unique minimal normal subgroup of $G_{uv}$,
namely $\soc(G_{uv})=\soc(G_v)$.
Since $G_{uv}\cong G_{vw}$,
$N:=\soc(G_{vw})\cong T^k$ is the unique minimal normal subgroup of $G_{vw}$.
Clearly, $N\cap \soc(G_v)\lhd G_{vw}$,
by the minimality of $N$, either $N\cap \soc(G_v)=1$ or $N=\soc(G_v)$.
For the former case,
we have
$$N=N/(N\cap \soc(G_v))\cong \soc(G_v)N/\soc(G_v)\le G_v/\soc(G_v)\le\Out(T)\times\S_k.$$
By Lemma~\ref{p-part},
we obtain $|T|_r^k=|N|_r<(k!)_r<r^{k\over r-1}$,
a contradiction.
Therefore, $\soc(G_{vw})=\soc(G_v)=\soc(G_{uv})$.
It follows
$$\soc(G_v)^g=\soc(G_{uv})^g=\soc(G_{uv}^g)=\soc(G_{vw})=\soc(G_v),$$
so $\soc(G_v)$ is normal in $\l G_v,g\r=G$,
that is, $G_v$ contains a nontrivial normal subgroup of $G$, it is a contradiction.\qed

To treat the candidates in Table 1, we first prove two lemmas.

\begin{lemma}\label{s=2}
Suppose $\Ga$ is $(G,2)$-arc-transitive.
Then for each prime $r$, we have
$$|T|_r<r|\phi_1(G_{uv}\cap M)|_r^2.$$
\end{lemma}

\proof Suppose $|\phi_1(G_{uv}\cap M)|_r=r^l$
and $G\cong\A_n.o$ with $o\le\ZZ_2$.
Then $G_v\cong T^k.(\Out(T)\times(\A_k.o))$
and $G_{uv}/(G_{uv}\cap M)\cong\A_k.o$,
so
$$|G_{uv}|_r\le |G_{uv}\cap M|_r|\A_k.o|_r\le |\phi_1(G_{uv}\cap M)|_r^k|\A_k.o|_r=r^{kl}({k!\over 2})_r|o|_r.$$
Since $\Ga$ is $(G,2)$-arc-transitive,
$G_v=G_{uv}G_{vw}$,
we obtain
$$|T|_r^k|\Out(T)|_r({k!\over 2})_r|o|_r=|G_v|_r\le |G_{uv}|_r^2\le r^{2kl}({k!\over 2})_r^2|o|_r^2.$$
By Lemma~\ref{p-part}, it follows $|T|_r^k\le r^{2kl}(k!)_r<r^{2kl+{k\over r-1}}$ ,
hence $|T|_r<r^{2l+{1\over r-1}}<r^{2l+1}$,
the lemma follows.\qed

\begin{lemma}\label{s=3}
Suppose $\Ga$ is $(G,3)$-arc-transitive. Then for each prime $r$,
we have
$$|T|_r^{2k}<r^{k\over r-1}|\phi_1(G_{uv}\cap M)|_r^{3k}|\Out(T)|_r.$$
\end{lemma}

\proof Suppose $|\phi_1(G_{uv}\cap M)|_r=r^l$
and $G\cong\A_n.o$ with $o\le\ZZ_2$.
Then by the proof of Lemma~\ref{s=2},
$|G_{uv}|_r\le r^{kl}({k!\over 2})_r|o|_r$,
and since $\val(\Ga)=|G_v:G_{uv}|$, we obtain
$$\val(\Ga)_r={|G_v|_r\over |G_{uv}|_r}
\ge {|T|_r^k|\Out(T)|_r({k!\over 2})_r|o|_r\over r^{kl}({k!\over 2})_r|o|_r}\ge{|T|_r^k\over r^{kl}}.$$
Since $\Ga$ is $(G,3)$-arc-transitive, by Lemma~\ref{s-ArcT},
$\val(\Ga)^3\mid |G_v|$.
It follows $|T|_r^{3k}$ divides $r^{3kl}|T|_r^k|\Out(T)|_r(k!)_r$,
then the lemma follows by Lemma~\ref{p-part}.\qed

We now analyse the candidates in Table 1 by the following two lemmas.
The proofs need many information of the orders and the outer automorphism groups
of simple groups, for those we refer to \cite[P. 18--20]{LPS90}.

\begin{lemma}\label{Cand-1}
Suppose $G_v$ satisfies Row $l$ of Table $1$, where $l\in\{7,9,11,12,14-18,21,22,$
$26-29\}$.
Then $s=1$.
\end{lemma}

\proof Suppose on the contrary, $s\ge 2$. We divided the proof into two cases.

\vskip0.1in
\noindent{\un{Row 7.}}{\hspace{5pt}} Then $T=\PSL_2(p)$ with $p=2^m-1$ a Mesenna prime,
and $\phi_1(G_{uv}\cap M)\le\ZZ_p:\ZZ_{p-1}$.
It follows $|T|_2=2^m$ and $|\phi_1(G_{uv}\cap M)|_2\le 2$.
By Lemma~\ref{s=2}, we obtain $2^m<2^3$, so $m\le 2$
and $T$ is soluble, a contradiction.

\vskip0.1in
\noindent{\un{Remaining rows.}} Then the simple groups $T$ are specific with no parameter,
and either $\phi_1(G_{uv}\cap M)\cap T$ or $\phi_1(G_{uv}\cap M)$ is
given in Table 1.
Since $\phi_1(G_{uv}\cap M)/(\phi_1(G_{uv}\cap M)\cap T)\le \Out(T)$,
we have $|\phi_1(G_{uv}\cap M)|_r\le |\Out(T)|_r|\phi_1(G_{uv}\cap M)\cap T|_r$.
Then a direct computation shows that the triple $(|T|_r,|\Out(T)|,|\phi_1(G_{uv}\cap M)|_r)$
for some prime $r$ lies in Table 2 (we remark that in Row 16 of Table 1,
$\phi_1(G_{uv}\cap M)\le \M_{22}:\ZZ_2$
by Atlas \cite{Atlas}, so $|\phi_1(G_{uv}\cap M)|_3=3^2$).
For each row in Table 2,
we always have $|T|_r\ge r|\phi_1(G_{uv}\cap M)|_r^2$,
by Lemma~\ref{s=2}, it is a contradiction.\qed

\begin{table}
\[\begin{array}{llllll} \hline
l & T & |T|_r &  |\Out(T)|  &  |\phi_1(G_{uv}\cap M)|_r & r  \\ \hline
9 & \PSL_3(3) & 3^3 & 2 & 3  & 3 \\
11 & \PSU_3(3) & 3^3 & 2 & 3 & 3 \\
12 & \PSU_3(5) & 5^3 & 6 & 5 & 5 \\
14 & \PSU_4(3) & 3^6 & 8 & 3^2 & 3 \\
15 & \PSU_5(2) & 3^5 & 2 & 3 & 3 \\
16 & \PSU_6(2) & 3^6 & 6 & 3^2 & 3 \\
17 & \PSp_4(7) & 7^4 & 2 & 7 & 7 \\
18 & \Sp_4(8) & 3^4 & 6 & 3 & 3 \\
21 & G_2(3) & 3^6 & 3 & \le 3^2 & 3 \\
22 & {^2F_4(2)'} & 3^3 & 2 & 3 & 3 \\
26 & HS & 5^3 & 2 & 5 & 5 \\
27 & McL & 3^6 & 2 & 3^2 & 3 \\
28 & Co_2 & 3^6 & 1 & 3^2 & 3 \\
29 & Co_3 & 3^7 & 1 & 3^2 & 3 \\
\hline
\end{array}\]
\caption{Triples $(|T|_r,|\Out(T)|,|\phi_1(G_{uv}\cap M)|_r)$ of `sporadic' cases in Lemma~\ref{Cand-1}.}
\end{table}

\begin{lemma}\label{Cand-2}
Suppose $G_v$ satisfies Row $l$ of Table $1$ with $l\in\{{1-6},8,10,13,19,20,23-25\}$.
Then $s\le 2$.
\end{lemma}

\proof Suppose on the contrary, $s\ge 3$.
We investigate each row in Lemma~\ref{Cand-2} in the following.

\vskip0.1in
\noindent{\un{Row 1.}}{\hspace{5pt}}
Then $T=\A_c$ and $\A_l\lhd \phi_i(G_{uv}\cap M)\le \S_l\times\S_{c-l}$,
where $5\le l<c$, $c$ is not a prime and $l$ is greater than or equal to the largest prime less than $c$.

If $c=6$, then $l=5$, so $|T|_3=3^2$, $|\Out(T)|_3=1$
and $|\phi_i(G_{uv}\cap M)|_3=3$.
It follows $|T|_3^{2k}=3^{4k}>3^{k/2}|\Out(T)|_3|\phi_i(G_{uv}\cap M)|_3^{3k}$,
contradicting Lemma~\ref{s=3}.

Thus assume $c>6$ in the following.
By Lemma~\ref{part-d}, we may suppose ${\ov{G_{uv}}}\le\S_k$ is transitive,
and hence $\phi_{i}(G_{uv}\cap M)=\phi_j(G_{uv}\cap M)$ for all $1\le i,j\le k$
by Lemma~\ref{part-d}.
For each $(a_1,\dots,a_k)\in M$,
since $G_v=G_{vw}G_{uv}$, we have
$(a_1,\dots,a_k)=(b_1,\dots,b_k)\s(c_1,\dots,c_k)\t$,
where $(b_1,\dots,b_k)\in G_{vw}\cap M$,
$(c_1,\dots,c_k)\in G_{uv}\cap M$,
and $\s,\t\in\S_k$.
It follows $(a_1,\dots,a_k)=(b_1c_{1^\s},\dots,b_kc_{k^{\s}})\s\t$,
and so $\s\t=1$ and $a_i=b_ic_{i^{\s}}$.
Consequently, $\S_c=\Aut(\A_c)\cong\phi_i(M)=\phi_i(G_{vw}\cap M)\phi_i(G_{uv}\cap M)$
because $c_{i^{\s}}\in \phi_{i^{\s}}(G_{uv}\cap M)=\phi_i(G_{uv}\cap M)$.
Now by \cite[P. 9, Corollary 5]{LPS90},
one of the following holds:
\begin{itemize}
\item[(i)] $\phi_i(G_{vw}\cap M)\supseteq T$.
\item[(ii)] $1\le c-l\le 5$ and $\phi_i(G_{vw}\cap M)$ is $(c-l)$-homogeneous on $c$ points.
\end{itemize}

\noindent Moreover, since $\A_l\lhd \phi_i(G_{uv}\cap M)\le \S_l\times\S_{c-l}$
and ${\ov{G_{uv}}}\le\S_k$ is transitive,
we obtain that either the minimal normal subgroups of $G_{uv}$ (so of $G_{vw}$ as $G_{uv}\cong G_{vw}$) are isomorphic to
$\A_l^h$ if $\phi_i(G_{uv}\cap M)\le \S_l$,
or isomorphic to $\A_l^h$ and $S^h$ if $\phi_i(G_{uv}\cap M)$ is not contained in $\S_l$,
where $h\le k$ and $S\le\S_{c-l}$.

For case (i), $G_{vw}$ has a minimal normal simple group isomorphic to $\A_c^t$ for some $t\le k$,
which is a contradiction.

Consider case (ii). Assume first $\phi_i(G_{vw}\cap M)$ is almost simple for some $i$,
say with socle $Q$.
Then $G_{vw}$ has a minimal normal subgroup isomorphic to
$Q^m$ for some $m\le k$,
so is $G_{uv}$. By the above assertion regarding the minimal normal subgroups of $G_{uv}$,
we obtain that either $Q\cong\A_5$ and $c-l=5$ or $Q\cong\A_l$.
For the former case, $\phi_i(G_{vw}\cap M)\le\S_5$ is $5$-homogeneous on $c$ points,
a contradiction.
For the latter case, if $n-l=1$, then
$\phi_i(G_{uv}\cap M)$ and $\phi_i(G_{vw}\cap M)$
are almost simple with the socles both isomorphic to $\A_l$,
one easily deduces from Lemma~\ref{AS-Factori}
that $c=6$ and $l=5$, also a contradiction.
If $n-l\ge 2$, as $\phi_i(G_{vw}\cap M)$ is almost simple with socle
$\A_l$ and is $(n-l)$-homogeneous on $c$ points with $c>l$,
by the classification of 2-homogeneous permutation groups which are not 2-transitive
(refer to \cite[Section 6.2]{Li-book}),
we obtain that $\phi_i(G_{vw}\cap M)$ is $2$-transitive. 
It then follows from \cite[Theorem 5.3(S)]{Cameron}
that $(c,l)=(6,5)$ or $(15,7)$.
The former contradicts to the assumption $c>6$,
and the latter contradicts to $c-l\le 5$,
yielding a contradiction.

Assume now $\phi_i(G_{vw}\cap M)$ is not almost simple for each $i=1,2,\dots,k$.
If $c-l=1$, then $\A_l\le\phi_i(G_{uv}\cap M)\le\S_{l}$
and the minimal normal subgroups of $G_{uv}$ and $G_{vw}$
are isomorphic to $\A_l^h$ with $h\le k$,
it leads to that there is $i$ such that $\A_l\le \phi_i(G_{vw}\cap M)\le \S_c=\S_{l+1}$,
which is impossible as $\phi_i(G_{vw}\cap M)$ is not almost simple.
Suppose $c-l\ge 2$.
Then $\phi_i(G_{vw}\cap M)$ is $2$-homogeneous on $c$ points,
since it is not almost simple,
we conclude that $\phi_i(G_{vw}\cap M)$ is of affine type,
namely $\soc(\phi_i(G_{vw}\cap M))\cong\ZZ_p^d$ with $p$ a prime
and $c=p^d$.
It follows that $G_{vw}$ and so $G_{uv}$
has a minimal normal subgroup isomorphic to $(\ZZ_p^{d})^m$ for some positive integer $m$.
However, by the above assertion of the minimal normal subgroups of $G_{uv}$,
we have $\ZZ_p^d\le\S_{c-l}$.
As $c-l\le 5$, we have $c=p^d\le 5$,
again a contradiction.

\vskip0.1in
\noindent{\un{Row 2.}}{\hspace{5pt}} Since $\PSL_2(5)\cong\A_5$,
the above discussion with $c=6$ draws a contradiction.

\vskip0.1in
\noindent{\un{Row 3.}} Then $|T|=|\PSp_{2m}(q)|=q^{m^2}\prod_{i=1}^m(q^{2i}-1)$
with $m,q$ even.
Set $q=2^e$.
Assume $em\ne 6$.
By Lemma~\ref{Zsigmondy},
$2^{em}-1=q^{m}-1$ has a primitive prime divisor, say $r$,
with $r>em$. Suppose $|q^{m}-1|_r=r^l$.
Notice that $r$ does not divide $q^{i}-1$
and $q^{m+i}-1$ (as $q^{m+i}-1=q^m(q^i-1)+(q^m-1)$) for each $1\le i\le m-1$.
we have
$|T|_r=|(q^m-1)(q^{2m}-1)|_r=r^{2l}$.
Since $\Omega_{2m}^-(q)\lhd\phi_1(G_{uv}\cap M)$,
$|\Omega_{2m}^-(q)|={1\over 2}q^{m(m-1)}(q^m+1)\prod_{i=1}^{m-1}(q^{2i}-1)$
and $r>em$, we conclude $|\phi_1(G_{uv}\cap M)|_r=|q^m-1|_r=r^l$.
As $q$ is even, $|\Out(T)|=e$ if $m\ge 3$,
and $|\Out(T)|=2e$ if $m=2$ (see \cite[P.18]{LPS90}),
so $r$ does not divide $|\Out(T)|$.
It then follows from Lemma~\ref{s=3} that
$r^{4kl}<r^{3kl+{k\over r-1}}$,
a contradiction.

Assume now $em=6$. Since $m$ is even, $(q,m)=(2,6)$ or $(2^3,2)$,
and $T=\PSp_{12}(2)$ or $\PSp_4(8)$ respectively.
For both cases, we have $|T|_7=7^2$, $|\phi_1(G_{uv}\cap M)|_7=7$,
and $|\Out(T)|=1$ if $m=6$ and $|\Out(T)|=6$ if $m=2$.
By Lemma~\ref{s=3}, we obtain
$7^{4k}<7^{3k+{k\over 6}}$,
also a contradiction.

\vskip0.1in
\noindent{\un{Row 4.}} Then $|T|=|\POmega_{2m+1}(q)|={1\over 2}q^{m^2}\prod_{i=1}^m(q^{2i}-1)$,
with $m$ even and $q=p^e$ an odd prime power.

If $(p,e,m)=(2^t-1,1,2)$ for some $t$, then one easily deduces
that $|T|_p=p^4$, $|\phi_1(G_{uv}\cap M)|_p=p^2$.
Since $\Out(T)=\ZZ_2$ in the case, by Lemma~\ref{s=3}, we have $p^{8k}<p^{6k+{k\over p-1}}$,
a contradiction.

Suppose $(p,e,m)\ne(2^t-1,1,2)$.
By Lemma~\ref{Zsigmondy},
$p^{em}-1=q^{m}-1$ has a primitive prime divisor $r>em$.
Then with similarly discussion as in the proof of {\un{Row 3}},
we have $|T|_r=(q^m-1)_r^2$,
$|\phi_1(G_{uv}\cap M)|_r=(q^m-1)_r$
and $(r,|\Out(T)|)=(r,2e)=1$.
It follows from Lemma~\ref{s=3}
that $(q^m-1)_r^{2k}<(q^m-1)_r^kr^{k\over r-1}$,
also a contradiction.

\vskip0.1in
\noindent{\un{Row 5.}} Then $|T|=|\POmega_{2m}^+(q)|={1\over (4,q^m-1)} q^{m(m-1)}(q^m-1)\prod_{i=1}^{m-1}(q^{2i}-1)$,
with $m$ even.
Suppose $q=p^e$ with $p$ a prime.

Assume first $(p,em)\ne (2,6)$.
By Lemma~\ref{Zsigmondy},
$p^{em}-1=q^{m}-1$ has a primitive prime divisor $r>em$.
Set $|q^{m}-1|_r=r^l$.
Since $q^{m+i}-1=q^m(q^i-1)+(q^m-1)$,
we have $r$ does not divide $q^{i}-1$
and $q^{m+i}-1$ for each $1\le i\le m-1$.
Then, as $m$ is even, it turns out
$|T|_r=|(q^m-1)^2|_r=r^{2l}$.
Since $|\Ome_{2m-1}(q)|={1\over 2}q^{(m-1)^2}\prod_{i=1}^{m-1}(q^{2i}-1)$
and $r>em$,
we have $|\phi_1(G_{uv}\cap M)|_r=(q^m-1)_r=r^l$.
Notice that $|\Out(T)|$ divides $24e$,
$(r,|\Out(T)|)=1$, then Lemma~\ref{s=3} implies $r^{4kl}<r^{3kl+{k\over r-1}}$,
a contradiction.

Now consider the case $(p,em)=(2,6)$.
Since $m$ is even,
we have $(m,q)=(2,8)$ or $(6,2)$,
and so $T=\POmega_4^+(8)$ or $\POmega_{12}^+(2)$
respectively. In particular, $|T|_7=7^2$
and $|\Out(T)|_7=1$ for both cases.
Since $\Omega_{2m-1}(q)\lhd \phi_1(G_{uv}\cap M)$
and $(m,q)=(2,8)$ or $(6,2)$,
one easily computes out $|\phi_1(G_{uv}\cap M)|_7=7$.
Then Lemma~\ref{s=3} leads to
$7^{4k}<7^{3k+{k\over 6}}$,
also a contradiction.

\vskip0.1in
\noindent{\un{Row 6.}} Then $|T|=|\PSp_4(q)|={1\over (2,q-1)}q^4(q^2-1)(q^4-1)$
and $|\Out(T)|=2e$, where $q=p^e$ is a prime power.
Since $\PSp_2(q^2)\lhd \phi_1(G_{uv}\cap M)$,
we have $\phi_1(G_{uv}\cap M)\cap T\le \PSp_2(q^2).\ZZ_2$. 
It follows $|\phi_1(G_{uv}\cap M)|$ divides
$2|\PSp_2(q^2)||\Out(T)|$ and hence divides $4eq^2(q^4-1)$.

If $(p,e)=(2^t-1,1)$, then $|T|_p=p^4$, $\Out(T)=\ZZ_2$
and $|\phi_1(G_{uv}\cap M)|_p$ divides $p^2$,
by Lemma~\ref{s=3}, we have $p^{8k}<p^{6k+{k\over p-1}}$,
a contradiction.
Similarly, if $(p,e)=(2,3)$, then $|T|_7=7^2$, $|\Out(T)|=6$
and $|\phi_1(G_{uv}\cap M)|_7=7$,
Lemma~\ref{s=3} implies
$7^{4k}<7^{3k+{k\over 6}}$, also a contradiction.

Thus assume now $(p,e)\ne (2^s-1,1)$ and $(2,3)$.
By Lemma~\ref{Zsigmondy}, $p^{2e}-1$ has a primitive prime divisor $r>2e$.
Set $(p^{2e}-1)_r=r^l$.
Then $|T|_r=r^{2l}$, $(r,|\Out(T)|)=1$ and $|\phi_1(G_{uv}\cap M)|_r\le r^l$.
It then follows from Lemma~\ref{s=3} that $r^{4kl}<r^{3kl+{k\over r-1}}$,
again a contradiction.

\vskip0.1in
\noindent{\un{Remaining rows.}}
With similar discussions as in the last paragraph of the proof of Lemma~\ref{Cand-1},
it is routine to compute out that the
triple $(|T|_r,|\Out(T)|,|\phi_1(G_{uv}\cap M)|_r)$ with $r$ a prime
is listed in Table 3.
For each row there, we always have $|T|_r^{2k}>r^{k\over r-1}|\phi_1(G_{uv}\cap M)|_r^{3k}|\Out(T)|_r$,
by Lemma~\ref{s=3}, it is a contradiction.\qed

\begin{table}
\[\begin{array}{llllll} \hline
l & T & |T|_r &  |\Out(T)|  &  |\phi_1(G_{uv}\cap M)|_r & r  \\ \hline
8 & \PSL_2(8) & 3^2 & 3 & 3  & 3 \\
10 & \PSL_6(2) & 3^4 & 2 & 3^2 & 3 \\
13 & \PSU_4(2) & 3^4 & 2 & 3^2 & 3 \\
19 & \Sp_6(2) & 3^4 & 1 & 3^2 & 3 \\
20 & \POmega_8^+(2) & 2^{12} & 6 & \le 2^7 & 2 \\
23 & \M_{11} & 3^2 & 1 & 3 & 3 \\
24 & \M_{12} & 3^3 & 2 & \le 3^2 & 3 \\
25 & \M_{24} & 3^3 & 1 & 3^2 & 3 \\
\hline
\end{array}\]
\caption{Triples $(|T|_r,|\Out(T)|,|\phi_1(G_{uv}\cap M)|_r)$ of `sporadic' cases in Lemma~\ref{Cand-2}.}
\end{table}

Summarize Lemmas~\ref{part-d}, \ref{lem-1}, \ref{Cand-1}
and \ref{Cand-2}, Lemma~\ref{part(d)} follows.

\section{Subgroups $(b), (e)$ and $(f)$}

Suppose Hypothesis~\ref{Hypo} holds. In this section, we treat the cases
where $G_{v}$ satisfies parts (b), (e) and (f),
and complete the proof of Theorem~\ref{Thm-1}.

We first consider parts (b) and (e).
For both cases, $G_v\cong(\S_m\wr\S_k)\cap G$,
where $n=mk$, $m>1$ and $k>1$ for part (b),
and $n=m^k$, $m\ge 5$ and $k>1$ for part (e).
Let $M\cong\S_m^k$ be the base group of $\S_m\wr\S_k$.

\begin{lemma}\label{part(b,e)}
Suppose $G_{v}$ satisfies parts $(b)$ or $(e)$
and $s\ge 2$. Then one of the following statements holds,
where $\phi_i(G_{\a\b}\cap M)$ denotes the projection of
$G_{\a\b}\cap M$ on the $i$-component of $M$,
and $(\a,\b)=(u,v)$ or $(v,w)$.

\begin{itemize}
\item[(i)] $\A_m\lhd \phi_i(G_{\a\b}\cap M)$ for $i=1,2,\dots,k$.
\item[(ii)] $\A_l\le\phi_i(G_{\a\b}\cap M)\le\S_l\times\S_{m-l}$,
where $m$ is not a prime, $p\le l<m$ with $p$ the maximal prime less than $m$,
and $i=1,2,\dots,k$.
\end{itemize}
\end{lemma}

\proof
Set $\ov{G_v}=G_vM/M$,
$\ov{G_{uv}}=G_{uv}M/M$
and $\overline{G_{vw}}=G_{vw}M/M$.
Since $s\ge 2$, $G_v=G_{uv}G_{vw}$,
so $\overline{G_v}={\overline{G_{uv}}}\hspace{1pt}{\overline{G_{vw}}}$.
By \cite[Lemma 2.3]{GLX17},
$\ov{G_{\a\b}}$
is a transitive subgroup of $\S_k$,
where $(\a,\b)=(u,v)$ or $(v,w)$.
Then Lemma~\ref{Wreath-Factor} implies
$\phi_1(G_{\a\b}\cap M)=\cdots=\phi_k(G_{\a\b}\cap M)$,
and $\pi(\A_m)\subseteq\pi(\phi_i(G_{\a\b}\cap M))$ for each $i=1,2,\dots,k$.

If $m\ge 5$, by \cite[Theorem 4(ii)]{LPS2000},
$\phi_i(G_{\a\b}\cap M)$ satisfies part (i) or (ii).
If $2\le m\le 4$,
since $\pi(\A_m)\subseteq\pi(\phi_i(G_{\a\b}\cap M))$
and $\phi_i(G_{\a\b}\cap M)\le\S_m$,
it is easy to see that
$\phi_i(G_{\a\b}\cap M)\rhd \A_m$ for $m=2$ and $3$,
and $\phi_i(G_{\a\b}\cap M)\rhd\A_4$ or equals to $\S_3$ for $m=4$,
that is, $\phi_i(G_{\a\b}\cap M)$ also satisfies part (i) or (ii).\qed

For a set $\Ome$, if $\Del_1,\dots,\Del_k$ are its subsets with equal size $m$,
such that $\Ome$ is the disjoint union of $\Del_1,\dots,\Del_k$,
then we call $(\Del_1,\dots,\Del_k)$ is
a $m$-{\it homogeneous partition} of $\Ome$.

\begin{lemma}\label{part-b}
Suppose $G_v$ satisfies part $(b)$ of Theorem~$\ref{Max-SubG}$.
Then $s\le 2$.
\end{lemma}

\proof
Then $G_v\cong(\S_m\wr\S_k)\cap G$ with $n=mk$, $m>1$ and $k>1$.
Suppose for a contradiction that $s\ge 3$.
By Lemma~\ref{part(b,e)}, without loss of generality, we may assume $\phi_i(G_{vw}\cap M)$
satisfies part (i) or (ii) of Lemma~\ref{part(b,e)}.

First suppose case (i) occurs.
Then $\A_m\lhd \pi_i(G_{vw}\cap M)$.
Notice that the action of $G$ on $V\Ga$
is permutation equivalent to the action of $G$ on
the set of all $m$-homogeneous partitions of $\Ome$,
we may set $v=(\Del_1,\dots,\Del_k)$,
a $m$-homogeneous partition of $\Ome$.
Then $u=v^{g^{-1}}=(\Del_1^{g^{-1}},\dots,\Del_k^{g^{-1}})$.
Since $u\ne v$,
without loss of generality, we may suppose
$\Del_1\ne\Del_1^{g^{-1}}$,
so $q:=|\Del_1\cap\Del_1^{g^{-1}}|<m$.
Since $\pi_1(G_{uv}\cap M)\le\S_m$ fixes both
$\Del_1^{g^{-1}}$ and $\Del_1$,
we have $\pi_1(G_{uv}\cap M)\le\S_q\times\S_{m-q}$,
which contradicts $\A_m\lhd \pi_1(G_{uv}\cap M)$.

Suppose now case (ii) occurs.
In this case, with almost the same discussions
as in the proof of Case {\un{Row 1}} of Lemma~\ref{Cand-2},
one may have a contradiction.\qed

\begin{lemma}\label{part(e)}
Suppose $G_v$ satisfies part $(e)$ of Theorem~$\ref{Max-SubG}$.
Then $s\le 2$.
\end{lemma}

\proof By assumption, $G_v\cong(\S_m\wr\S_k)\cap G$, with $n=m^k$, $m\ge 5$ and $k>1$.
Suppose on the contrary that $s\ge 3$.
By Lemma~\ref{part(b,e)}, we may assume $\phi_i(G_{vw}\cap M)$
satisfies part (i) or (ii) of Lemma~\ref{part(b,e)}.

For case (i), since $m\ge 5$, $\A_m$ is nonabelian simple.
Notice that in this case the arguments in the proof of Lemma~\ref{lem-1}
are efficient, then one may have that $G_v$ has a subgroup $\soc(G_v)\cong\A_m^k$
which is normal in $G$,
a contradiction.
For the case (ii), with similar discussions as in the proof of case {\un{Row 1}}
of Lemma~\ref{Cand-2},
one may has a contradiction.\qed

We finally consider part $(f)$.

\begin{lemma}\label{part(f)}
Suppose $G_v$ satisfies part $(f)$ of Theorem~$\ref{Max-SubG}$.
Then $(\soc(G_v),\soc(G_{uv}))=(P\Omega_8^+(q),\Omega_7(q))$
with $q$ an odd prime power and
$s\le 2$.
\end{lemma}

\proof By assumption, $G_v$ is almost simple,
then by \cite[Corollary 3.4]{GLX17},
we have $s\le 2$.
Set $T=\soc(G_v)$.
Since $G_v=G_{uv}G_{uv}^g$,
by Lemma~\ref{AS-Factori},
$(T,\soc(G_{uv}))=(\A_6,\A_5)$,
$(\M_{12},\M_{11})$, $(\Sp_4(2^f),\Sp_2(4^f))$
with $f\ge 2$,
or $(P\Omega_8^+(q),\Omega_7(q))$.

Suppose $T\cong\A_6$.
Since $G_v$ is maximal in $G$, we have $n=7$ or $8$.
If $n=7$, then $(G,G_v)=(\A_7,\A_6)$ or $(\S_7,\S_6)$,
so $G$ is 2-transitive on $V\Ga$,
and hence $\Ga\cong\K_7$ is an undirected complete graph, a contradiction.
If $n=8$, then $(G,G_v)=(\A_8,\S_6)$,
and the action of $G$ on $V\Ga$ is permutation equivalent to the action
of $G$ on $\Omega^{\{2\}}$, the set of $2$-subsets of $\Ome$.
It easily follows (or see \cite[TABLE B.2]{DM})
that $G$ is of rank 3 and the two nontrivial orbits
are with length $12$ and $15$, so $\Ga$ is undirected, a contradiction.

Suppose $G_v\cong\M_{12}$. Since $G_v$ is maximal in $G$,
we have $n=12$, and $G\cong\A_{12}$ is 2-transitive on $V\Ga$,
hence $\Ga=\K_{12}$ is an undirected complete graph,
a contradiction.

If $T\cong\Sp_4(2^f)$, by part (II)(B) of \cite[T{\tiny HEOREM}]{LPS87},
there is a subgroup $H\cong\S_m\wr\S_2$ with $m={1\over 2}q^2(q^2-1)$
such that $G_v<H<G$, so $G$ is not primitive on $V\Ga$, a contradiction.
Similarly, if $T\cong P\Omega_8^+(q)$ with $q$ even, 
there is a subgroup $H$ with $\soc(H)\cong\Sp_8(q)$
such that $G_v<H<G$ by part (II)(B) of \cite[T{\tiny HEOREM}]{LPS87}, hence $G$ is not primitive on $V\Ga$,
again a contradiction.\qed

Now, we may complete the proof of Theorem~\ref{Thm-1}.

\vskip0.1in
\noindent{\bf Proof of Theorem~\ref{Thm-1}.}
Since the full automorphism groups of the directed cycles are soluble,
$\val(\Ga)\ne 1$. If $\val(\Ga)=2$, by \cite[Theorem 5]{Neumann},
$\Aut\Ga$ is a dihedral group, a contradiction.
Thus $\Ga$ is of valency at least three.
Since $G_v$ with $v\in V\Ga$ is a maximal subgroup of $G$,
$G_v$ satisfies one of parts (a)-(f) of Theorem~\ref{Max-SubG}.
If $G_v$ satisfies parts (a) and (c), by Lemmas~\ref{part-a}
and \ref{part-c},
we have $s=1$. If $G_v$ satisfies parts (b), (d), (e) and (f), by Lemmas~\ref{part-b},
\ref{part(d)},
\ref{part(e)} and \ref{part(f)}, we have $s\le 2$.
This completes the proof of Theorem~\ref{Thm-1}.\qed


\begin{thebibliography}{99}


\bibitem{Magma}
W. Bosma, J. Cannon, C. Playoust, The MAGMA algebra system I:
The user language, {\it J. Symbolic Comput.} \textbf{24} (1997), 235--265.

\bibitem{BB82}
N. Blackburn, B. Huppert, Finite Groups $II$, Springer-verlag, New York, 1982.



\bibitem{Cameron}
P. J. Cameron, Finite permutation groups and finite simple groups,
{\it Bull. London Math. Soc.} {\bf 13} (1981), 1-22.

\bibitem{CPW93}
P. J. Cameron, C. E. Praefer, N. C. Wormald,
Infinite highly arc-transitive digraphs and universal covering digraphs,
{\it Combinarorica} {\bf 13} (1993), no. 4, 377--396.

\bibitem{CLP95}
M. Conder, P. Lorimer, C. E. Praeger,
Constructions for arc-transitive digraphs,
{\it J. Austral. Math. Soc. Ser. A} {\bf 59} (1995), no. 1, 61--80.


\bibitem{Atlas}
J. H. Conway, R. T. Curtis, S. P. Noton, R. A. Parker and R. A. Wilson,
{\it Atlas of Finite Groups}, Clarendon Press, Oxford, 1985.

\bibitem{DM}
\J. D. Dixon, B. Mortimer,
Permutation groups, Springer-Verlag, New York, 1997.

\bibitem{Evans97}
D. M. Evans, An infinite highly arc-transitive digraphs,
{\it Europ. J. Combin.}
{\bf 18} (1997), no. 3, 281--286.


\bibitem{GLX18}
M. Giudici, C. H. Li, B. Z. Xia,
An infinite family of vertex-primitive $2$-arc-transitive digraphs,
{\it J. Combin. Theory Ser. B} {\bf 17} (2017), 1--13.


\bibitem{GLX17}
M. Giudici, C. H. Li, B. Z. Xia, Vertex-primitive $s$-arc-transitive digraphs of linear groups,
arXiv:1710.09518v1.


\bibitem{MX18}
M. Giudici, B. Xia, Vertex-quasiprimitive 2-arc-transitive digraphs,
{\it Ars Math. Contemp.} {\bf 14} (2018), no. 1, 67--82.

\bibitem{Li-book}
C. H. Li, Permutation groups and symmetrical graphs,
The University of Western Australia, 2014.


\bibitem{LPS87}
M. W. Liebeck, C. E. Praeger, J. Saxl,
A Classification of the Maximal Subgroups of the
Finite Alternating and Symmetric Groups,
{\it J. Algebra} {\bf 111} (1987), 365--383.


\bibitem{LPS90}
M. W. Liebeck, C. E. Praeger, J. Saxl,
The maximal facorizations of the finite simple groups and their automorphism groups,
{\it Mem. Amer. Math. Soc.}
{\bf 86} (1990), Volume 86, Number 432.


\bibitem{LPS2000}
M. W. Liebeck, C. E. Praeger, J. Saxl,
Transitive subgroups of primitive permutation groups,
{\it J. Algebra} {\bf 234} (2000), 291--361.

\bibitem{MMSZ02}
A. Malni\v c, D. Maru\v si\v c, N. Seifter, B. Zgrabli\' c,
Highly arc-transitive digraphs with no homomorphism onto $Z$,
{\it Combinatorica} {\bf 22} (2002), 435--443.

\bibitem{Neumann}
P. M. Neumann, Finite permutation groups, edge-coloured graphs and matrices,
{\it Topics in group theroy and computation}
(Proc. Summer School,
University College, Galway, 1973),
pp. 82--118, Academic Press, London, 1977.

\bibitem{Praeger89}
C. E. Praeger, Highly arc-transitive digraphs, {\it Europ.
J. Combin.} {\bf 10} (1989), no. 3, 281--292.

\bibitem{Praeger90}
C. E. Praeger, Finite primitive permutation groups:
a survey, {\it Groups-Canberra 1989}, 63--84, Lecture Notes in Math. 1456, Springer, Berlin, 1990.


\bibitem{Weiss}
R. Weiss, The nonexistence of $8$-transitive graphs,
{\it Combinatorica} {\bf 1} (1981), no. 3, 309--311.
\end{thebibliography}
\end{document}